\documentclass[12pt]{amsart}
\usepackage{amsfonts}  
\usepackage{mathrsfs}
\usepackage{pictex, color}

\theoremstyle{plain}
\newtheorem{Thm}{Theorem}[section]
 
\newtheorem{Lem}{Lemma}[section]

\newtheorem{Cor}{Corollary}[section]
\theoremstyle{definition}

\newtheorem{Rmk}{Remark}

\newcommand{\ve}{\varepsilon}

 \newcommand{\Q}{\mathbb{Q}}

\def\eps{{\varepsilon}}

\title{Transcendence with Rosen continued fractions}
\author{Yann Bugeaud}
\address{Universit\'e de Strasbourg\\
Math\'ematiques\\
 7, rue Ren\'e Descartes\\
 67084 Strasbourg cedex,  France}
\email{bugeaud@math.unistra.fr}

\author{Pascal Hubert}
\address{LATP, case cour A, Facult\'e des sciences Saint J\'er\^ome \\
Avenue Escadrille Normandie Niemen \\
13397 Marseille cedex 20, France}
\email{hubert@cmi.univ-mrs.fr}
\author{Thomas A. Schmidt}
\address{Department of Mathematics\\Oregon State University\\Corvallis, OR 97331, USA}
\email{toms@math.orst.edu}
\thanks{The second named author is partially supported by project blanc ANR: ANR-06-BLAN-0038.   The third author thanks FRUMAN, Marseille and the Universit\'e  P. C\'ezanne.} 
\date{24 June 2010}

\usepackage{hyperref}

\begin{document}
\baselineskip=17pt

\begin{abstract}
We give the first transcendence results 
for the  Rosen continued fractions.   
Introduced over half a century ago, these fractions  
expand real numbers in terms of certain algebraic numbers.
\end{abstract}

\maketitle

\section{Introduction} 
In 1954, D. ~Rosen defined 
an infinite family of continued fraction algorithms \cite{R}.     
Introduced to aid in the study of certain Fuchsian groups,   
these continued fractions were applied some thirty  
years later by J.~Lehner \cite{Le} in the study of 
Diophantine approximation by orbits of these groups. 

The Rosen continued fractions and variants   have been of recent interest, leading to results especially about their dynamical and arithmetical properties, see \cite{BKS}, \cite{Nzero}, \cite{DKS}; as well on their applications to the study of geodesics on related hyperbolic surfaces, see \cite{SS}, \cite{BS}, \cite{MS};   and to Teichm\"uller geodesics arising from (Veech) translation surfaces, see \cite{Vch}, \cite{AH}, \cite{SU}  and \cite{AS}.   
Several basic questions remain open, including that of arithmetically characterizing the real numbers having a finite Rosen continued fraction expansion, 
see \cite{Leu}, \cite{TetAl} and  \cite{AS}.
Background on Rosen continued fractions is given in the next section.

The first transcendence criteria for regular continued fractions
were proved by E. Maillet, H. Davenport and K. F. Roth, A. Baker, and recently
improved by B. Adamczewski and Y. Bugeaud, see
\cite{AB05,AB07} and the references given there.
In particular, Theorem 4.1 of \cite{AB07} asserts that, if $\xi$ is an algebraic irrational
number with sequence of convergents $(p_n/q_n)_{n\ge 1}$, then the
sequence $(q_n)_{n \ge 1}$ cannot increase too rapidly.
It is natural to ask whether similar
transcendence results can be proven using  Rosen continued fractions.   
We give the first such results.

\begin{Thm}\label{t:transc}   
Fix $ \lambda = 2 \cos \pi/m$ for an   
integer $m >3\,$,  and denote the field extension degree $[\mathbb Q(\lambda): \mathbb Q]$ by $D$.    
If a real number $\xi \notin  \mathbb Q(\lambda)$ has an infinite  
expansion  in Rosen continued fraction
over $\mathbb Q(\lambda)$ of  convergents $p_n/q_n$ satisfying
$$
\limsup_{n \to \infty}\,  \dfrac{\log\log q_n}{n} > \log (2 D - 1)\,,
$$ 
then $\xi$ is transcendental.  
\end{Thm}

For stating our second result, we associate to the Rosen continued fraction expansion
\[
[\, \ve_1(x):r_1(x),\,
\ve_2(x):r_2(x),\ldots ,\, \ve_n(x):r_n(x),\ldots ] : = \dfrac{\varepsilon_1}{r_1 \lambda + \dfrac{\varepsilon_2}{r_2 \lambda + \cdots}}   
\]
of a real number  $x$
the  sequence of pairs of integers  $(\eps_i, r_i)_{i\ge 1}$, 
which we call the {\em partial quotients}, 
and thus consider ${\mathcal A} = \{\pm 1\} \times \mathbb N$ as the alphabet of 
the Rosen continued fraction expansions.

As usual,  we denote the length of a finite word $U = u_1 \cdots u_k$ as $|U| = k$.
  For any positive integer $s$, we write
$U^s$ for the word $U \ldots U$ ($s$ times repeated concatenation
of the word $U$). More generally, for any positive real number
$s$, we denote by $U^s$ the word
$U^{\lfloor s \rfloor}U'$, where $U'$ is the prefix of
$U$ of length $\left\lceil(s - \lfloor s \rfloor)\vert U\vert\right\rceil$.

Just as Adamczewski and  Bugeaud \cite{AB05} showed 
for regular continued fraction expansions,  a real number  whose Rosen 
continued fraction expansion is appropriately ``stammering''  must be transcendental.

\begin{Thm}\label{t:wordGrowth}
Fix $ \lambda = 2 \cos \pi/m$ for an 
integer $m >3\,$,  
and denote the field extension degree $[\mathbb Q(\lambda): \mathbb Q]$ by $D$.
Let $\xi$ be an infinite Rosen continued fraction 
with convergents $(p_n/q_n)_{n \ge 1}$ such that
\[
B := \limsup_n q_n^{1/n} < + \infty\,.
\]
Write 
\[
b := \liminf_n q_n^{1/n}\,.
\]
Assume that there are two infinite
sequences $(U_n)_{n \ge 1}$ and $(V_n)_{n \ge 1}$
of finite words over the alphabet ${\mathcal A}$ 
and an infinite sequence $(w_n)_{n \ge 1}$
of real numbers greater than $1$ such that, for $n \ge 1$,
the word $U_n V_n^{w_n}$ is a prefix of
the infinite word composed of the partial quotients
of $\xi$. If
\begin{equation}\label{eq:limSupBound}
\limsup_{n \to + \infty} \, 
{|U_n| + w_n |V_n| \over 2 |U_n| + |V_n|} 
> {3 D \over 2} \cdot { \log B \over \log b}\,, 
\end{equation}
then $\xi$ is either (at most)
quadratic over ${\mathbb Q} (\lambda)$ or is transcendental.
\end{Thm}

\noindent
Lemma ~\ref{l:bPos}  gives that $\log b$ is positive.\\

The key to our proofs is that both numerator and denominator 
of a Rosen convergent dominate their respective conjugates 
in an appropriate fashion, see Lemma ~\ref {l:domConjugates}.  
From this one can bound the height of a Rosen convergent 
in terms of its denominator, see Lemma~\ref{l:heightBd}.      
Then, exactly as in the case of regular continued fractions, we apply
tools from Diophantine approximation, namely an extension to number fields
of the Roth theorem, for the proof of Theorem \ref{t:transc},
and the Schmidt Subspace Theorem for the proof of Theorem \ref{t:wordGrowth}.

Both theorems  are weaker than their analogues for regular 
continued fractions, since we have to work in a number field of degree $D$ rather than
in the field $\mathbb Q$.
However, for $m=4$ and $m=6$, that is, for $\lambda_4 = \sqrt{2}$ and    
$\lambda_6 = \sqrt{3}$, our results can be considerably strengthened and we can get,  
essentially, the exact analogues of the results established for the regular  
continued fractions. 
The key point is that, in both cases, for every  convergent $p_n/q_n$, exactly one of $p_n, q_n$ is in $ \mathbb Z$,  the other being in $\lambda \mathbb Z\,$, see Remark 2 below.

\section{Background}
\subsection{Rosen fractions} We set $\lambda = \lambda_m = 2 \cos \frac{\pi}{m}$ and
${\mathbb I}_m = [-\lambda/2, \lambda/2\,)$ for $m
\ge 3$. For a fixed integer $m \ge 3$, the Rosen continued fraction
map is defined by
\[
T(x) = \begin{cases} \left| \frac{1}{x} \right| \, - \,
\lambda \lfloor \left| \,\frac{1}{\lambda x} \right| + \frac{1}{2}
\rfloor & x \ne 0; \\
\\
  0   &  x = 0
\end{cases}
\]
for $x \in {\mathbb I}_m$; here and below, we  omit the index
``$m$" whenever it is clear from context.  For $n \geq 1$, we define 
\[
\ve_n(x) = \ve (T^{n-1}x) \qquad \mbox{and} \qquad  
r_n(x) = r(T^{n-1}x)
\]
with
\[
\ve (y) = \mbox{sgn}(y) \qquad \mbox{and} \qquad  
r(y) = \left\lfloor \, \left| \frac{1}{\lambda y} \right| + \frac{1}{2} 
\right\rfloor .
\]
Then, as Rosen showed in \cite{R},  the Rosen continued fraction expansion of $x$ 
is given by
\[
[\, \ve_1(x):r_1(x),\,
\ve_2(x):r_2(x),\ldots ,\, \ve_n(x):r_n(x),\ldots ] : = \dfrac{\varepsilon_1}{r_1 \lambda + \dfrac{\varepsilon_2}{r_2 \lambda + \cdots}}\,.\]
  As usual we define the {\em convergents} $p_n/q_n$
of $x\in {\mathbb I}_m$ by
\[
\begin{pmatrix}
p_{-1}  &  p_0  \\
q_{-1}  &  q_0
\end{pmatrix}
\, = \,
\begin{pmatrix}
1  &  0  \\
0  &  1
\end{pmatrix}
\]
and
\[
\begin{pmatrix}
p_{n-1}  &  p_{n}  \\
q_{n-1}  &  q_{n}
\end{pmatrix}
\, = \,
\begin{pmatrix}
0 &  \ve_{1}  \\
1  &  \lambda r_{1}
\end{pmatrix}
\begin{pmatrix}
0 &  \ve_{2}  \\
1  &  \lambda r_{2}
\end{pmatrix}
\cdots
\begin{pmatrix}
0 &  \ve_{n}  \\
1  &  \lambda r_{n}
\end{pmatrix}
\]
for $n \ge 1$. From this definition it is immediate that $\left|
p_{n-1}q_{n} \, - \, q_{n-1}p_n \right| =1$, and that the
well-known recurrence relations
\begin{eqnarray*}
p_{-1}=1; \, p_0=0; & p_n=\lambda r_n p_{n-1}+\varepsilon_n p_{n-2},\,
n\geq 1\\
q_{-1}=0; \, q_0=1; & q_n=\lambda r_n q_{n-1}+\varepsilon_n q_{n-2},\,
n\geq 1,
\end{eqnarray*}
hold. It also follows that

\begin{equation}
\begin{pmatrix}
p_{n-1}  &  q_{n-1}  \\
p_{n}  &  q_{n}
\end{pmatrix}
\, = \,
\begin{pmatrix}
0 &  1  \\
\ve_{n}  &  \lambda r_{n}
\end{pmatrix}
\begin{pmatrix}
0 &  1  \\
\ve_{n-1}  &  \lambda r_{n-1}
\end{pmatrix}
\cdots
\begin{pmatrix}
0 &  1  \\
\ve_{1}  &  \lambda r_{1}
\end{pmatrix} ,
\end{equation}
giving
\begin{equation}
\frac{p_n}{q_n} \, = \,[\, \ve_1:r_1,\,
\ve_2:r_2,\ldots ,\, \ve_n:r_n ]
\end{equation}
and
\begin{equation}\label{eq:mirror}
\frac{q_{n-1}}{q_n} \, = \,[\, 1:r_n,\,
\ve_n:r_{n-1},\ldots ,\, \ve_2:r_1 ]\; .
\end{equation}

We define 
\begin{equation}\label{eq:matrix} M_n = \begin{pmatrix}
p_{n-1}  &  p_{n}  \\
q_{n-1}  &  q_{n}
\end{pmatrix}, \, \text{and find that}\;  x = M_n \cdot T^n(x),  
\end{equation}
where $\cdot$ denotes the usual fractional linear operation, namely  
$$
x = \frac{p_{n-1} T^n(x) + p_n}{q_{n-1} T^n(x) + q_n}. 
$$

\subsection{Approximation with Rosen fractions}
We briefly discuss the convergence of the ``convergents'' to $x$.     One can rephrase some of Rosen's original arguments in terms of  the (standard number theoretic) natural extension map $\mathcal T(x,y) = (\, T(x), \frac{1}{r \lambda + \varepsilon y}\,)\,$ where $r = r_1(x)$ and $\varepsilon = \varepsilon_1(x)\,$.   The ``mirror formula'' Equation ~(\ref{eq:mirror})  shows that  $\mathcal T^n(x,y) = (\, T^n(x),  \frac{q_{n-1}}{q_n}\,)\,$.  Extending earlier work of H. ~Nakada, 
it is shown in \cite{BKS}   
that $\mathcal T(x,y)$ has planar domain with $y$-coordinates between $0$ and $R = R(\lambda)\,$,  where $R= 1$ if the index $m$ is even  and, otherwise, 
$R$ is the positive root of $R^2 + (2 - \lambda) R - 1 = 0\,$,  
in which case we have $1 > R > \lambda/2\,$ (see Lemma 3.3 of \cite{BKS}).    
Therefore, the sequence $(q_n)_{n \ge 1}$ is strictly increasing.     
But,   as Rosen mentions,   if $x$ has infinite expansion,  then either $\ve_n = 1$ or $r_n >1$ occurs infinitely often;  from this one has both that $q_n\ge1$ for all $n$ and that the limit as $n$
tends to infinity 
of $q_n$ is infinite.\\    

One easily adapts Rosen's arguments so as to find the following.  

\begin{Lem} \label{l:bPos}
For every $x \in \mathbb I_m$ of infinite expansion, we have
\[ 
\liminf_n q_n^{1/n} > 1\,.\]  
\end{Lem}

\begin{proof}
We know that the sequence $(q_n)_{n \ge 1}$ increases and that, 
if either $\eps_n = 1$ or $r_n >1$,  then   $q_n > \lambda q_{n-1} $.  Furthermore, there are no more than $h$ consecutive indices $i$
with $(\eps_i, r_i) = (- 1,1)$, with $h = m/2, (m-3)/2$ depending on the parity of $m$ (see \cite{R} or \cite{BKS}).
Consequently, for any $n$,   there is some $i = 1, \ldots , h+1$  such that $q_{n+i} > \lambda q_{n + i-1} $, giving 
$$
q_{n+h+1} \ge q_{n+i} > 
\lambda q_{n+i-1} \ge \lambda q_n\,.
$$
As $q_1 \ge \lambda$, letting $s(n) = 1+ \bigg\lfloor\dfrac{n-1}{h+1}\bigg\rfloor$, we have   $q_n \ge \lambda^{s(n)}$.   Since $\lambda > 1$, 
this proves the lemma. 
\end{proof}

\begin{Rmk}  In fact,  H.~Nakada \cite{N} shows that {\em for almost all} such $x$,  $ \lim_{n\to \infty} \dfrac{1}{n} \log q_n$ 
exists, being  equal to  one half of 
the  entropy of $T$.   He also shows that the entropy equals $C \cdot (m-2) \pi^2/(2 m)$, where $C =  1/\log(1+R)$  when $m$ is odd, and equals $1/\log[(1 + \cos \pi/m)/\sin \pi/m]$ when $m$ is even. 
This $C$ is the normalizing constant to give a probability measure on the domain of the planar natural extension $\mathcal T$, see \cite{BKS}.  
\end{Rmk} 

\bigskip

    Rosen also gave  bounds on $\vert\, x - {p_n}/{q_n}\,\vert\,$.    
    Using Equation~\eqref{eq:matrix} (as in Nakada \cite{N})
\[\left\vert\, x - \dfrac{p_n}{q_n}\,\right\vert\, = \dfrac{1}{q_{n}^{2}} \, 
\dfrac{1}{|\, \frac{q_{n+1}}{q_n} + T^{n+1} x\,|}\,,\]   
one has the easy lower bound 
\[ \dfrac{1}{q_{n}\, (  q_{n+1} + q_n)} < \left\vert\, x - \dfrac{p_n}{q_n}\,\right\vert\,\,.\]   
One also finds in this manner
\begin{equation}\label{eq:upperBdApprox} \left\vert\, x - \dfrac{p_n}{q_n}\,\right\vert 
< \dfrac{c_1}{q_n q_{n+1}}\,,
\end{equation}
with $c_1 = c_1(\lambda) = 2/(2 - R \lambda)\,$.       Thus convergence does hold.\\

Theorems 4.4 and 4.5 (depending on parity of $m$) of \cite{BKS} give  $c_2$ such that 
$\vert\, x - {p_n}/{q_n}\,\vert\,  < c_2/q_{n}^{2}\,$,  with the upper bound 
$c_2 \le 1/2+ \lceil m/4\rceil$.

\subsection{Traces in Hecke groups}
Rosen introduced his continued fractions to study the Hecke  groups.   The Hecke (triangle Fuchsian) group $G_m$ with $m \in \{3, 4, 5, \dots\,  \}$ is the group generated by 
\[
 \begin{pmatrix}
	1   &  \lambda_m\\
 	0   &  1
\end{pmatrix}\, \, \text{and} \;\; 
 \begin{pmatrix}
        0   &  -1\\
        1   &  0
\end{pmatrix},\]
with $\lambda_m$ as above.  The Rosen expansion of a real number terminates at a finite term if and only if $x$ is a
parabolic fixed point of  $G_m\,$, see \cite{R}.   These points are clearly contained in $\mathbb Q(\lambda_m)\,$ but in general there are elements of this field that have infinite Rosen expansion,  see \cite{Leu}, \cite{TetAl} and  \cite{AS}.

\begin{Rmk}   The values of finite Rosen expansion form the set $G_{q} \cdot\infty$, which is in fact a subset of  $\lambda \, \mathbb Q(\lambda^{2}) \cup
\{\infty\}$.   To see this, one uses  induction on word length in the  generators displayed above --- an ordered pair $(a, c)$ giving a column of any element of $G_q$ must be such that exactly one element of the pair is in $\mathbb Z[\lambda^2]\,$, and the other is  in $\lambda \mathbb Z[\lambda^2]\,$.       Note that this also applies to  convergents $p_n/q_n$:  exactly one of $p_n, q_n$ is in $ \mathbb Z[\lambda^2]$, the other being in $\lambda \mathbb Z[\lambda^2]\,$. 
\end{Rmk} 

\bigskip 
 
When $q=3\,$, we have    $G_3 = \text{PSL}(2, \mathbb Z)\,$.    In general each $G_m$ is isomorphic to the free product of a cyclic group of order two   and a cyclic group of order $m$.            Recall that a Fuchsian triangle group is generated by {\em even} words in the reflections about the sides of some hyperbolic triangle.   Thus any Fuchsian triangle group is of index two in the group generated by these reflections;  for each $G_m$,  we denote  this larger group by $\Delta_m\,$.

Since $\lambda_m$ is the sum of the root of unity  $\zeta_{2 m} := \exp{2 \pi i/(2 m)}$ with its complex conjugate,   $\mathbb Q(\lambda_m)$ is a number field of degree $d := \phi( 2 m)/2$ over the rationals, where $\phi$ denotes the Euler totient function.    
 
The following key phenomenon property of Hecke groups can be shown in various manners.     The result holds for a larger class of groups, from  Corollary 5 of \cite{SW},  due to \cite{CW} (extending the arguments  from  $G_m$ to $\Delta_m$ is straightforward).   Independent of this earlier work,  Bogomolny-Schmit \cite{BS} gave a clever proof of the result specifically for $\Delta_m\,$.   See the next remark for another perspective.

\begin{Thm}\label{t:bigTrace}    Fix $m$ as above, and let $\Delta_m$ be the full reflection group in which $G_m$ has index two.    Then for any $M \in \Delta_m$ whose trace is of absolute value greater than $2$, we have 
\[ |\, \text{tr}(M)\,| \ge  |\, \sigma(\,\text{tr}(M)\,)\,|\,,\]
where $\sigma$ is any field embedding of $\mathbb Q(\lambda_m)\,$.
\end{Thm}

\begin{Rmk}   This result can be proven ``geometrically''.    Up to conjugacy,  each of the Hecke groups appears as the Veech group of some translation surface, 
see \cite{Vch};  the elements whose trace is of absolute value at least 2 are the ``derivatives'' of the affine pseudo-Anosov diffeomorphisms of the surface.     The dilatation of a pseudo-Anosov $\phi$  is the dominant eigenvalue $\lambda$ of the action of $\phi$ on the integral homology of the underlying surface.   (The other eigenvalues are  hence conjugates of $\lambda$.)   The corresponding element of the Veech group has  trace of absolute value $\lambda + \lambda^{-1}$ from which it follows that this trace dominates its conjugates.  
\end{Rmk}

\subsection{Approximation by algebraic numbers}

The following result was announced by Roth \cite{Roth}
and proven by LeVeque, see Chapter 4 of \cite{L}.  
(The version below is Theorem 2.5  of \cite{B}.)     
Recall that given an algebraic number $\alpha$, its {\em naive height}, denoted 
by $H(\alpha)\,$, is the largest  absolute value  of the coefficients of 
its minimal polynomial over $\mathbb Z\,$.

\begin{Thm}\label{t:Lev}(LeVeque) Let $K$ be a number field,  and $\xi$ a real algebraic number not in $K\,$.  Then, for any $\epsilon >0\,$, there exists a positive constant $c(\xi, K, \epsilon)$ such that 
\[ |\, \xi - \alpha \,| > \dfrac{c(\xi, K, \epsilon)}{H(\alpha)^{2+\epsilon}} \]
holds for every $\alpha$ in $K$.
\end{Thm}

 \bigskip         
     
          The {\em logarithmic Weil height} of $\alpha$ lying in a number field $K$  of degree $D$ over $\mathbb Q$ is  $h(\alpha) =  \frac{1}{D}\sum_{\nu} \, 
\log^{+} \max_{\nu \in M_K}\{|| \alpha ||_{\nu}\}$, where $\log^{+} t$ equals $0$ if $t\le 1$ and $M_K$ denotes the places (finite and infinite ``primes'') of the field,  and $|| \cdot ||_{\nu}$ is the $\nu$-absolute value.    
This definition is independent of the field $K$ containing $\alpha$. 
The two heights are related by 
\begin{equation}\label{eq:height}    
\log H(\alpha) \le \deg(\alpha) h(\alpha) + \log 2, 
\end{equation}
for any non-zero algebraic number $\alpha$,
see Lemma 3.11 from \cite{W}. \\

We recall a consequence of the W.~ Schmidt Subspace Theorem.

\begin{Thm}\label{t:subspace}   
Let $d$ be a positive integer and $\xi$ be a real algebraic
number of degree greater than $d$. Then, for every positive $\eps$,
there exist only finitely many algebraic numbers $\alpha$
of degree at most $d$ such that
$$
|\xi - \alpha| < H(\alpha)^{-d-1-\eps}.
$$
\end{Thm}
\noindent
Note that the Roth theorem is exactly the case $d=1$ of   
Theorem ~\ref{t:subspace}.

In the proof of Theorem ~\ref{t:wordGrowth}, we could apply Theorem ~\ref{t:subspace}, but  
the algebraic numbers $\alpha$ which we use to
approximate $\xi$   are of degree
at most $2$  over a fixed number field. In this situation,
the next theorem, kindly communicated to us by Evertse \cite{E}, yields  
a stronger result than the previous one.  

\begin{Thm}(Evertse)\label{t:evertse}    
Let $K$ be a real algebraic number field of degree $d$.
Let $t$ be a positive integer and $\xi$ be a real algebraic
number of degree greater than $t$ over $K$. 
Then, for every positive $\eps$,
there exist only finitely many algebraic numbers $\alpha$
of degree $t$ over $K$ and $\delta$ over $\Q$ such that
$$
|\xi - \alpha| < H (\alpha)^{-dt(t+1+\eps)/\delta}.
$$
\end{Thm}

 Note that Theorem ~\ref{t:evertse} extends Theorem ~\ref{t:Lev}.

\subsection{Sturmian sequences: towards an application of Theorem 1.2} 

To give an explicit  family of Rosen expansions satisfying the hypotheses of Theorem ~1.2, we recall a result of  \cite{AB10} on Sturmian sequences.

 Let $a$ and $b$ be letters in some alphabet.  The complexity function of a sequence ${\bf u} = u_1u_2 \cdots$ with values in $\{a,b\}$ is given by letting  $p(n,{\bf u})$ be the number of distinct words of length $n$ that occur in ${\bf u}$.   A sequence ${\bf u}$ is called {\em Sturmian} if its 
 complexity satisfies $p(n,{\bf u}) = n+1$ for all $n$.     As Arnoux \cite{A} writes, one can obtain any such sequence by taking a ray with irrational slope in the real plane and intersecting it with an integral grid,  assigning $a$ when the ray intersects a horizontal grid line and $b$ when it meets a vertical grid line.    Indeed, the {\em slope} of a Sturmian sequence is the density of $a$  in the sequence  (one shows that the limit as $n$ tends to infinity  
of the average of the number of occurrences $a$ in $u_1 \cdots u_n$ exists, see \cite{A}, Proposition 6.1.10).

\begin{Lem}\label{l:sturm}
Let ${\bf u}$ be a Sturmian word whose slope has an unbounded
regular continued fraction expansion. Then, for every positive integer $n$, there
are finite words $U$, $V$ and a positive real number $s$ such that
$U V^s$ is a prefix of ${\bf u}$ and $|U V^s| \ge n |U V|$.
\end{Lem}

\begin{proof}
This follows from the proof of Proposition 11.1
from \cite{AB10}. 
\end{proof} 

\begin{Rmk}   
We apply the above lemma to Sturmian sequences where both $a,b$ are of the form $(\eps , r)$, with $\eps = \pm 1$ and $r \in \mathbb N$.   In particular,  we use this in the context of  Rosen expansions to prove Corollary \ref{c:sturm}.   
\end{Rmk}

\section{Bounding the height of convergents} 
 
In what follows, we fix $\lambda= \lambda_m$ for some $m > 3$, and suppose that $\xi \in (0, \lambda/2)\,$ is a real algebraic number having an infinite  Rosen continued fraction expansion over $\mathbb Q(\lambda)\,$.   Our goal is to estimate the naive height $H( {p_n}/{q_n}\,)\,$ 
of the $n$th convergent $p_n/q_n$.     
In light of Theorem ~\ref{t:bigTrace}, we let $n_0$ be the  least value of $n$ such that $q_n > 2\,$.      

\begin{Lem}\label{l:domConjugates}   Let $c_3 = c_3(\lambda)$ be defined by $c_3 = \min_{\sigma}  \frac{|\, \sigma(\lambda)\,|}{\lambda}\,$, where the minimum is taken over all field embeddings of $\mathbb Q(\lambda)$ into $\mathbb R\,$.     Then for all  $n\ge n_0$,  and any such $\sigma\,$, 
 we have both 
\[ q_n \ge c_3 \; \vert\,\sigma(q_n)\,\vert \;\; \text{and}\;\; p_n \ge c_3 \; \vert\,\sigma(p_n)\,\vert\,.\]
\end{Lem}
\begin{proof}    For any $n\ge n_0$, recall that $M_n = \begin{pmatrix}
p_{n-1}  &  p_n  \\
q_{n-1}  &  q_n
\end{pmatrix}\,$;   this is clearly an element of $\Delta_m\,$.     By  Theorem ~\ref{t:bigTrace} we have $q_n + p_{n-1}   \ge \vert\, \sigma(q_n + p_{n-1} )\,\vert\,$.

Now let $j \in \mathbb N\,$ and set
\[ M_{n,j} = \begin{pmatrix}
p_{n-1}  &  p_n  \\
q_{n-1}  &  q_n
\end{pmatrix}     \begin{pmatrix}
1  &  j \lambda  \\
0  &  1
\end{pmatrix}   =  \begin{pmatrix}
p_{n-1}  &  p_n + j \lambda p_{n-1}  \\
q_{n-1}  &  q_n+ j \lambda q_{n-1} 
\end{pmatrix}\,.\]
This is also an element of  $\Delta_m$ of trace greater than 2, and hence 
 \[ |\, p_{n-1} +  q_n+ j \lambda q_{n-1} \,|  \ge |\, \sigma( p_{n-1} +  q_n)+ j \sigma(\lambda q_{n-1}  )\,|\,.\]
Since this holds for {\em all} positive $j\,$, we must have that $\lambda q_{n-1} \ge \vert\, \sigma(\lambda q_{n-1}  )\,\vert\,$.  That is,  
\[q_{n-1} \ge \frac{|\, \sigma(\lambda)\,|}{\lambda} \,\vert \sigma(q_{n-1}  )\vert\,\ge \biggl( \min_{\sigma}  \frac{|\, \sigma(\lambda)\,|}{\lambda} \biggr)\, \vert\sigma(q_{n-1}  )\vert\,.\]

Similarly,  using
\[ N_{n,j} = \begin{pmatrix}
p_{n-1}  &  p_n  \\
q_{n-1}  &  q_n
\end{pmatrix}     \begin{pmatrix}
1  & 0 \\
 j \lambda  &  1
\end{pmatrix}   =  \begin{pmatrix}
p_{n-1}  + j \lambda p_{n}&  p_n  \\
q_{n-1} + j \lambda q_{n} &  q_n
\end{pmatrix}\,,\]
we find 
\[p_{n} \ge \frac{|\, \sigma(\lambda)\,|}{\lambda} \, \vert\sigma(p_{n}  )\vert\,\ge \biggl( \min_{\sigma}  \frac{|\, \sigma(\lambda)\,|}{\lambda} \biggr)\, \vert\sigma( p_{n}  )\vert\,.\]
\end{proof} 

\begin{Rmk}  We conjecture that in fact $q_n$ is always greater than or equal to its conjugates, thus that in the above one can replace $c_3$ by $1$.
\end{Rmk}

\begin{Lem}\label{l:heightBd}   Let $D$ denote the field extension degree  
$[\mathbb Q(\,\lambda\,):\mathbb Q\,]\,$.
There exists a constant $c_4 = c_4 (\lambda)$ such that for all $n \ge n_0$, 
\[
H(p_n/q_n) \le c_4 q_n^D.
\]
\end{Lem}
 
\begin{proof} 
Since $p_n$ and $q_n$ are algebraic integers
of degree at most $D$, it follows from Lemma ~\ref{l:domConjugates} that
\[
h(p_n/q_n) \le \sum_{\sigma} \, \frac{1}{D}
\log \max\{|\sigma(p_n)|, |\sigma(q_n)|\}  \le c'_4 + \log q_n,
\]
where $\sigma$ runs through the complex embeddings, for a suitable
positive constant $c'_4$.
Using ~\eqref{eq:height}, we get the asserted estimate. 
\end{proof} 

\begin{Lem}\label{l:htPeriodic} 
Let $\alpha$ be a real number in
$[-\lambda/ 2, \lambda / 2)$ with an ultimately periodic
expansion in Rosen continued fraction.
Denote by $(p_n / q_n)_{n \ge 1}$ the sequence of its convergents.
Denote by $\mu$ the length of the preperiod
and by $\nu$ the length of the period, with the
convention that $\mu=0$ if the expansion is purely periodic.
Then $\alpha$ is of degree
at most $2$ over $\mathbb Q(\lambda)$, and there exists $c_5 = c_5(\lambda, \alpha)$ such that
\[
H(\alpha) \le c_5 (q_{\mu} q_{\mu+\nu})^{D}.
\]
\end{Lem}

\begin{proof} In the notation of Equation ~\eqref{eq:matrix},    $\alpha$ is fixed by $M = M_{\mu}^{-1}M_{\mu + \nu}$. It thus satisfies a quadratic equation with entries in $\mathbb Z[\lambda]$,  and hence is of degree at most $2$ over $\mathbb Q(\lambda)$.   Indeed,   $\alpha$ is a root of $f(x) = cx + (d-a)x - b$ with $a,b,c,d$ denoting the entries of $M$.   Each entry is a $\mathbb Z$-linear combination of monomials of the form $r s$ with $r$ an entry of $M_{\mu}$ and $s$ an entry of $M_{\mu + \nu}$.       

 Now,   $\alpha$ is also a root of $\tilde f(x)  = \prod_{\sigma}\, \sigma(f) (x) \in \mathbb Z[x]$, where $\sigma(f)$ denotes the result of applying  $\sigma$ to the coefficients of $f(x)$.   By Lemma  ~\ref{l:domConjugates},  all of the conjugates of each of  $p_{\mu}, p_{\mu-1}, q_{\mu-1}$, $q_{\mu}$ can be bounded by the product of $q_{\mu}$ with a constant depending upon $\alpha$ and $\lambda$.  Similarly for the entries of $M_{\mu + \nu}$.  
After some computation,
we conclude that the height of $\alpha$ is  $\ll q_\mu^{D} q_{\mu+\nu}^D$.
(One checks that the case of $\mu=0$ is subsumed by the above.) 
\end{proof}

\begin{Rmk}  Whereas a real number whose regular continued fraction expansion is ultimately periodic  is exactly of degree two over the field of rational numbers,  in the previous lemma the words ``at most'' are necessary.   Indeed, 
$x = 1$ has  an ultimately periodic Rosen expansion with respect to any $\lambda_m$ with $m$ even, \cite{R}.  Further examples of  elements of $\mathbb Q(\lambda_m)$ with periodic expansions are easily given when $m \in \{4,6\}$,  see Corollary 1 of \cite{SS}.   Yet further examples, including cases with $m \in \{7, 9\}$,  are given in \cite{RT}, \cite{TetAl}.

\end{Rmk}

\section{Transcendence results}

As usual, $\ll$ and $\gg$ denote inequality with implied constant.

\subsection{Applying Roth--LeVeque:  the proof of Theorem ~\ref{t:transc}}  

We now show that the sequence of denominators of convergents to an algebraic number cannot grow too quickly.   Theorem ~\ref{t:transc} then follows.

\begin{proof}  Let $\eps$ be a positive real number.   Let $\zeta$ be an algebraic number
having an infinite Rosen expansion with convergents $r_n / s_n$.

By the Roth--LeVeque Theorem \ref{t:Lev},  we have 
\[
|\zeta - r_n / s_n| \gg  H(r_n / s_n)^{- 2 - \eps}, \quad \hbox{for $n \ge 1$}.
\]
And, hence by Lemma ~\ref{l:heightBd},  
for $n \ge n_0 = n_0(\zeta)$, we have $|\zeta - r_n / s_n| \gg s_n^{- 2D - D \eps}\,$.
Inequality ~\eqref{eq:upperBdApprox}  
then gives that there exists a constant $c_6$ (independent of $n\ge n_0$) such that 
\[
s_{n+1} < c_6 s_n^{2D - 1 + D \eps}\,.
\]

Set $a = 2D-1 + D \eps\,$.  For $j < n_0$, define $\ell_j$ such that $s_j < \ell_j s_{j-1}^a\,$.   
We set $c_7 = \max\{ 1, c_6,  \ell_1, \dots, \ell_{n_0-1}\}$ and find that for any $n>1$ 
\[ s_{n+1} < c_7 s_{n}^{a} < c_7(c_7 s_{n-1}^a)^a \le (c_7 s_{n-1})^{a^2}\]
and continuing in this manner,  we have $s_{n+1} <  (c_7 s_1)^{a^n}$.  Since $s_{n+1} > s_n$ letting $c_8 = c_7 s_1$,   gives  $ \log s_{n} < a^n \log c_8\,$.   From this follows that 
$$
\limsup_{n \to + \infty} \, {\log \log s_n \over n} < \log (D(2 + \eps) - 1).  
$$ 
Letting $\eps$ go to zero, we see that every algebraic number satisfies
$$
\limsup_{n \to + \infty} \, {\log \log s_n \over n} \le  \log (2 D - 1),
$$
as asserted.
\end{proof}

\bigskip 
\subsection{Proof of Theorem  ~\ref{t:wordGrowth} and an application.}

\begin{proof}  With $\lambda = 2 \cos \pi/m$ fixed, given  $\xi$ of infinite   
Rosen continued fraction 
with convergents $(p_n/q_n)_{n \ge 1}$, we let $b = \liminf_n q_n^{1/n}$ and $B = \limsup_n q_n^{1/n}$, and assume that $B <\infty$.   
Let $\eta$ be a positive real number with $b-1< \eta < b$.  
Since there are only finitely many $n$ with either $q_n^{1/n} < b-\eta$ or $q_n^{1/n} > B+ \eta$, we have both that  $q_n \gg (b-\eta)^n$ and $q_n \ll (B+\eta)^n$\,.

Suppose that  $w$ is a positive real number and  $U$, $V$ are finite words in $\{\pm 1\} \times \mathbb N$ such that   $UV^w$
is a prefix of
the infinite word composed of the partial quotients
of $\xi$. Denote by $\alpha$ the real number 
of degree at most two over $\Q (\lambda)$ whose 
Rosen continued fraction is given by the word $U V^{\infty}$, where $V^{\infty}$
means the concatenation of infinitely many copies of $V$.   
Set $|U|=u$ and $|V| = v$. 
Since $\xi$ and $\alpha$ have their first $\lfloor u + v w \rfloor$ 
partial quotients in common, we have
$$
|\xi - \alpha| < c_2  \,q_{\lfloor u + v w \rfloor}^{-2}
\ll (b-\eta)^{-2(u + v w)}.
$$
Furthermore, it follows from Lemma ~\ref{l:htPeriodic} that
\[
H(\alpha) \ll (q_{u}\,q_{u+v})^D \ll (B+ \eta)^{D (2u + v)}.
\]
Combined with the previous inequality, this gives
\[
|\xi - \alpha| \ll H(\alpha)^{- 2 (u + v w) \log (b - \eta) /
(D(2u + v) \log (B + \eta))}.
\]

Now suppose that $\xi$ is algebraic of degree greater than two over $\mathbb Q(\lambda)$.    
Then, for every $\eps > 0$,  
there exists a positive constant $C(\eps)$
such that every real algebraic number $\beta$ of degree at most $2$ over $\mathbb Q(\lambda)$
satisfies
\[
|\xi - \beta| > C(\eps) H(\beta)^{-3-\eps}.  
\]
This follows from Theorem \ref{t:Lev} if $\beta$ is in $Q(\lambda)$ and, otherwise, by applying
Theorem ~\ref{t:evertse}  with $t=2$ and $d t = \delta$ to each
subfield $K$ of $\mathbb Q(\lambda)$. 

This proves that $\xi$ must be transcendental if there are
$u, v, w$ such that $u + vw$ is arbitrarily large and
$$
{2 (u + v w) \log b \over D(2u + v) \log B} > 3,  
$$
as asserted. 
\end{proof}

\begin{Cor}\label{c:sturm}
A Rosen continued fraction whose sequence of partial quotients 
is Sturmian with slope of unbounded  regular continued fraction partial quotients represents a transcendental number.
\end{Cor}

\begin{proof}
Combine Lemma  ~\ref{l:sturm} with the previous Theorem. 
\end{proof}

\begin{Rmk}
Using the Subspace Theorem as in \cite{AB05} does not yield in general
an improvement of Theorem ~\ref{t:wordGrowth}. In case $u=0, b=B$, inequality \eqref{eq:limSupBound}
reduces to $w > 3D/2$, while, proceeding as in \cite{AB05},
we would get $w > 2D - 1$. However, if $b$ is much smaller
than $B$ and $D$ is small, then the approach of \cite{AB05} presumably
gives a slightly better result than Theorem ~\ref{t:wordGrowth}.
\end{Rmk}
 

\end{document}